\mag=\magstep1 
\input epsf 
\documentclass[10pt,a4paper]{amsart} 
\usepackage{amssymb,latexsym} 
\usepackage{verbatim} 

\def\Cal{\mathcal} 
 
\def\frak{\mathfrak}

\def\<<{\langle } 
\def\>>{\rangle }

\numberwithin{equation}{section} 
\newtheorem{theorem}{Theorem}[section] 
\newtheorem{proposition}[theorem]{Proposition} 
\newtheorem{corollary}[theorem]{Corollary} 
\newtheorem{definition}[theorem]{Definition}

\newtheorem{lemma}[theorem]{Lemma}

\newtheorem{notation}[theorem]{Notation}

\def\Coker{\operatorname{Coker}}

\def\<{\langle} 
\def\>{\rangle}

\begin{document} 
 
\title{Algebraic correspondences between genus three curves and
certain Calabi-Yau varieties} 
 
\author{Tomohide Terasoma} 
 
\address{Graduate School of Mathematical Sicences, 
the University of Tokyo, 
3-8-1 Komaba, Meguro, 
Tokyo, 153-8914, Tokyo, Japan}

\maketitle 
 
\makeatletter 
\renewcommand{\@evenhead}{\tiny \thepage \hfill 
Algebraic correspondences
\hfill}

\renewcommand{\@oddhead}{\tiny \hfill 
T. Terasoma
 \hfill \thepage}

\section{Introduction}
\label{intro section}

Let
$
\Phi_L:\bold P^3--\to \bold P(L^*)
$ 
be a net of quadrics, i.e. two dimensional linear system of quadrics.
If the linear system $\Phi_L$ is sufficiently generic, then
\begin{enumerate}
\item
the base locus
$B$
consists of eight points
$b_1, \dots, b_8$,
and 
\item
the locus $C$ of singular quadrics in $\bold P(L)$
is a smooth plane quartic curve.
\end{enumerate}

A set $B$ of eight points in $\bold P^3$ of this form is called
a regular Cayley octad and the plane quartic curve $C$ is called 
the associated curve of $B$. The moduli space of regular Cayley octads
is isomorphic to the moduli space of plane quartic curves with
markings of even theta divisors (cf. \cite{DO}).

Let $B=\{b_1, \dots, b_8\}$ be a regular Cayley octad,
$H_1, \dots, H_8$ hyperplanes
in the dual projective space 
$\widehat{\bold P^3}=\{(\xi_1:\xi_2:\xi_3:\xi_4)\}$
of $\bold P^3$
corresponding to the points $b_1, \dots, b_8$, and 
$X$ the double covering of $\widehat{\bold P^3}$
branched along the eight hyperplanes $H_1, \dots, H_8$.
In this paper, we study a certain algebraic correspondence 
between the associated curve $C$ and the double covering $X$ of 
$\widehat{\bold P^3}$ based on a geometry of curves
on $X$. Since the variety $X$ has only quotient singularities,
$H^i(X,\bold Q)$ is equipped with a pure Hodge structure
of weight $i$.
The first main theorem of this paper is the following.

\begin{theorem}[Theorem \ref{correspondence thm}]
\label{anounce1}
Under the above notations,
we have an injection of Hodge structures
\begin{equation}
\label{algebraic correspondece map}
cyl: H^1(C, \bold Q)(-1)\to H^3(X,\bold Q)
\end{equation}
induced by an algebraic 
correspondence.
In particular, the Hodge structure $H^3(X,\bold Q)$
is not irreducible.
\end{theorem}

The same algebraic correspondence is constructed independently
in unpublished work by I.Dolgachev-E.Merkmann (1994). 
The author would
like to express his thanks to 
Prof. I.Dolgachev for the communication.
We remark that the above homomorphism specializes to that defined by
the ``symmetric construction'' in \cite{T1},
if the curve $C$ specializes to a hyperelliptic curve. 
Let us recall the symmetric construction briefly.
Let $C$ be a hyperelliptic curve of genus three 
defined by 
$$
y^2=\prod_{i=1}^8(x-\lambda_i).
$$
We define a double covering $X$ of $\widehat{\bold P^3}$ by
$$
\eta^2=\prod_{i=1}^8(\lambda_i^3\xi_1+\lambda_i^2\xi_2+
\lambda_i\xi_3+\xi_4)
$$
in the weighted projective space $\bold P(1,1,1,1,4)$.
Here $\xi_0,\xi_1,\xi_2,\xi_3$ are homogeneous coordinates of
degree one and $\eta$ is that of degree four.
Then the variety $X$ is a Calabi-Yau variety, i.e. it 
has only Gorenstein singularities and has
the trivial dualizing sheaf admitting a global crepant 
resolution (see \cite{CM}).

Let $C^{(1)},C^{(2)},C^{(3)}$ be three copies of the curve $C$.
The action of $\mu_2=\langle \iota \rangle$ on $C$ defined by 
the hyperelliptic involution $\iota$
yields an action of ${\mu_2}^3$ on $C^3=C^{(1)}\times C^{(2)}\times C^{(3)}$.
The symmetric group $\frak S_3$ of degree three acts on 
$C^3$ by permuting the components, and as
a consequence, we obtain an action of semi-direct product 
$({\mu_2}^3)\rtimes \frak S_3$ on $C^3$. 
Let $N$ be the kernel of the homomorphism 
${\mu_2}^3 \to \mu_2:(\zeta_1,\zeta_2,\zeta_3)\mapsto \zeta_1\zeta_2\zeta_3$.
Then the semi-direct product 
$G=N\rtimes \frak S_3$ is a normal subgroup of ${\mu_2}^3\rtimes \frak S_3$
of index two.

Since the quotient of $C^3$ by
the group $({\mu_2}^3)\rtimes \frak S_3$ is isomorphic to 
$(\bold P^1)^3/\frak S_3\simeq \bold P^3$, the quotient of $C^3$
by $G$ is a double covering of $\bold P^3$.
Since the covering map
$$
C^3/G \to
C^3/({\mu_2}^3)\rtimes \frak S_3 \simeq \bold P^3
$$
branches exactly along the eight hyperplanes 
$H_i=\{\lambda_i^3\xi_1+\lambda_i^2\xi_2+
\lambda_i\xi_3+\xi_4=0\}$, 
we have an isomorphism
\begin{equation}
\label{basic isom}
X \overset{\simeq}\to C^3/G
\end{equation}
over $(\bold P^1)^3/\frak S_3\simeq \bold P^3$.
The above isomorphism is written as
\begin{align*}
&(\xi_1:\xi_2:\xi_3:\xi_4:\eta) \\
&=
(-1:x_1+x_2+x_3:-(x_1x_2+x_2x_3+x_3x_1):
x_1x_2x_3:y_1y_2y_3).
\end{align*}
on the affine part of $C^3$.
Here $(x_i,y_i)$ is a coordinate of
the affine part of $C^{(i)}$ for $i=1,2,3$.

Since the variety $X$ has only rational singularities,
$H^3(X,\bold Q)$ is equipped with a pure Hodge structure of weight three.
By the isomorphism (\ref{basic isom}),
we have an isomorphism of Hodge structures :
\begin{equation}
\label{exterior structure}
\bigwedge^3H^1(C, \bold Q)\simeq H^3(X, \bold Q).
\end{equation}
The composite morphism
$$
\begin{matrix}
cyl:&C^{(1)}\times C & \to & C^{(1)}\times C^{(2)}\times C^{(3)}
& \overset{\varphi}\to & X \\
&p_1\times p& \mapsto &
p_1\times \iota(p_1)\times p
\end{matrix}
$$
defines a family of smooth 
rational curves in $X$ parameterized by $C$.
(See Section \ref{construction of family}
for details.)
From this family of smooth rational curves in $X$, we get a homomorphism
of Hodge structures:
\begin{equation}
\label{cylinder for he}
cyl:H^1(C, \bold Q)(-1)\to H^3(X,\bold Q).
\end{equation}
This map is the specialization of the map 
(\ref{algebraic correspondece map}) in the following sense.
Let $M$ be the moduli space of genus three curves with level 
$2N$-structures ($N\geq 2$)  and
$p:\Cal C\to M$ and $\Xi:\Cal X \to M$ the universal families of
plane quartic curves and double 
coverings of $\bold P^3$. The family of
algebraic correspondences in Theorem \ref{anounce1} give rise to
a homomorphism of 
variations of Hodge structures:
\begin{equation}
\label{universal cokernel}
cyl:R^1p_*\bold Q(-1)\to R^3\Xi_*\bold Q.
\end{equation}
The variations of Hodge structures $R^1p_*\bold Q(-1)$ and $\bold
R^3\Xi_*\bold Q$
and the homomorphism $cyl$
extend to those on $M$.
The homomorphism (\ref{cylinder for he})
is the fiber of
(\ref{universal cokernel}) at the point corresponding to 
the hyperelliptic curve $C$.

In the last part of this paper, we study the variation 
of Hodge structures obtained by the cokernel $\Coker(cyl)$ of the
homomorphism (\ref{universal cokernel}).
By the isomorphism (\ref{exterior structure}), 
the restriction of $\Coker(cyl)$ to the hyperelliptic locus
is isomorphic to the primitive part of
the third higher direct image of the relative Jacobian scheme
as variations of Hodge structures.
The Hodge type of $\Coker(cyl)$ at each fiber is $(1,6,6,1)$. 
Nevertheless, we have
the following theorem.
\begin{theorem}[Theorem \ref{revised second main theorem}]
\label{anounce2}
There exists no polarized abelian scheme $a:\Cal A \to M$ of relative
 dimension three, whose degree three primitive part
$$
R^3a_{*,prim}\bold Q=\Coker(R^1a_*\bold Q(-1)\overset{L}\to R^3a_*\bold Q)
$$ 
of the higher direct image is isomorphic to $\Coker(cyl)$.
Here $L$ denotes the Lefschetz operator for the polarization.
\end{theorem}

The contents of the paper is as follows.
We recall several basic facts concerning net of quadrics in 
Section \ref{recall quadratic complex}.
In Section \ref{construction of correspondence},
we construct an algebraic correspondence 
between a non-hyperelliptic curve of genus three
and a double covering of $\widehat{\bold P^3}$ 
using the theory of net of quadrics.
In Section 4 and 5, we prove the injectivity of the homomorphism
induced by this correspondence
(Theorem \ref{anounce1}). In Section \ref{ivhs for double covering}
and Section \ref{computer calculation},
we compute infinitesimal variations of Hodge
structure to conclude Theorem \ref{anounce2}.

{\bf Acknowledgment}. 
The author would like to thank Keiji Matsumoto for discussions,
which gives a motivation of this paper.
The statement of Theorem \ref{anounce2}
is a refinement of the first version of this paper,
which is suggested by the referee.
The author would like to express his thanks to the referee.
After finishing this paper, the relevant paper \cite{GSZ} was
pointed out by K.Zuo and B.van Geemen, which treats general 
configurations of eight hyperplanes and hyperelliptic locus.
The author is grateful for the information.
\begin{notation}
\label{cor quad and symm}
Quadratic polynomials 
$Q(x_1, \dots, x_k)$ of 
$x_1, \dots, x_k$ are in one to one
correspondence with $k\times k$ symmetric
matrices $Q'$ by the relation
$Q(x_1, \dots, x_k)=(x_1, \dots, x_k){Q'}{}^t(x_1, \dots, x_k)$.
The corresponding symmetric matrix $Q'$
is also denoted by $Q$ if
there are no confusions. 
The rank of the quadric is defined by the rank of the corresponding
symmetric matrix. The space of quadratic form on $V$ is denoted by $Sym^2(V)$.

For a vector space $V$, the projective
space $(V-\{0\})/\bold C^{\times}$ associated to $V$ is denoted as $\bold P(V)$.
For a homogeneous polynomial $f$ on $V$, 
the subvariety of $\bold P(V)$ defined by $f$
is denoted as $Z(f)$.
The point of the dual projective space $\bold P(V^*)$
corresponding to a hyperplane 
$W\subset \bold P(V)$ is denoted by $[W]$.
The hyperplane in $\bold P(V^*)$ associated to 
a point $b\in \bold P(V)$
is denoted by $M_b$.

\end{notation}

\section{Net of quadrics in the three dimensional projective space}
\label{recall quadratic complex}

Let $V$ be a four dimensional vector space.
For a three dimensional subspace $L$ in $Sym^2(V)$, 
the associated linear system $\bold P(V) --\to \bold P(L^*)$
is denoted as $\Phi_L$. A two dimensional linear system is called a net.
The member of the linear system corresponding to 
$t\in \bold P(L)$ is denoted as
$Q_t\subset \bold P(V)$.
The locus of singular quadrics in $\bold P(L)$ is denoted as $C(L)$.
A net of quadrics in $\bold P(V)$ is called regular if
$C(L)$ is smooth. In this case,
$C(L)$ is a smooth plane quartic curve and 
is called the plane quartic curve associated to $L$.
By choosing a basis $Q_1,Q_2,Q_3$ of $L$
and basis $x_0, \dots, x_3$ of $V$,
the subvariety $C(L)$ of $\bold P(L)$ is defined by the polynomial
$\det(t_1Q_1+t_2Q_2+t_3Q_3)$. 
For a regular net of quadrics $L$,
the rank of the quadric $t_1Q_1+t_2Q_2+t_3Q_3$ is three for all 
$t=(t_1:t_2:t_3)\in C(L)$ and
the base locus $B$ of the net of quadrics $\Phi_L$ consists of
distinct eight points.
Moreover any four points in the base locus
do not lie on a common plane.

Let $b_1, \dots, b_7$ be generic seven points in $\bold P(V)$.
Then the space $L$ of quadratic polynomials on $V$ vanishing at
$b_1, \dots,b_7$ is a three dimensional vector space,
and $\Phi_L$ is a regular net of quadrics.
Moreover the base locus $B$ of the net of quadrics 
is a zero dimensional reduced subscheme consisting of eight points 
containing $\{b_1, \dots, b_7\}$.

A configuration of eight points in $\bold P(V)$ 
is called a regular 
Cayley octad if it is the base locus of a regular net of quadrics.
By associating the base locus a net of quadrics, 
we get a one to one 
correspondence between
the set of regular nets of quadrics in 
$\bold P(V)$ and
that of regular Cayley octads. 
The set of choices of two points in $\{b_1, \dots, b_8\}$
corresponds one to one to the set of odd theta divisors of $C(L)$, 
and there is one to one
correspondence between the set of odd theta divisors and that
of bitangents of $C(L)$.

We have the following four moduli spaces.
\begin{enumerate}
\item
The moduli space $M_1$ of regular Cayley octads 
$\{b_1, \dots, b_8\}$ in $\bold P(V)$. 
\item
The moduli space $M_2$ of regular nets of quadrics 
in $\bold P(V)$.
\item
The moduli space $M_3$ of
smooth plane quartic curves.
\item
The moduli space $M_4$ of non-hyperelliptic smooth curves of genus three.
\end{enumerate}
Then we have the following morphisms
$$
M_1\overset{f_1}\to M_2 \overset{f_2}\to M_3 \overset{f_3}\to M_4.
$$
The morphisms $f_1$ and $f_3$ are isomorphisms and 
the morphism $f_2$ is an etale finite morphism of degree 36.
\begin{definition}
[Steinerian curve]
Let $\<Q_1,Q_2,Q_3\>$ be a regular net of quadrics and
$C$ the associated plane quartic curve.
Then $Z(Q_t)$ is a cone over a smooth conic in $\bold P(V)$
for $t \in C$.
By attaching the vertex $s(t)$ of $Z(Q_t)$ to a point $t\in C$,
we have a map $s:C \to \bold P(V)$.
The image $s(C)$ of $C$ under the morphism $s$ is 
called the Steinerian curve (\cite{DO}). 
\end{definition}
It is known that the Steinerian curve is of degree 6 and
isomorphic to 
$C$ via the map
$s:C \to \bold P(V)$.
The Steinerian curve is disjoint from the base locus
of the net of quadrics.
The embedding $s$
is equal to the linear system defined by
the sum of canonical class and one of 36 even theta divisors.

\section{Double covering branched along eight hyperplanes}
\label{construction of correspondence}
We use the same notations as in the last section.
Let $\{b_1, \dots, b_8\}$ be a regular Cayley octad
in $\bold P(V)$
and $C$ be the associated plane curve.
Let 
$H_1, \dots, H_8$ be the hyperplanes corresponding to $b_1, \dots, b_8$ in 
the dual projective space $\bold P(V^*)$ of $\bold P(V)$.
Let $X$ be the double covering branched
along the union of the hyperplanes
$H_1, \dots, H_8$. The covering map
$X \to \bold P(V^*)$ is denoted by 
$\pi$.
Since any four points in $\{b_1, \dots, b_8\}$ are not contained in
a common plane by the condition of regular Cayley octad,
the union $D=\cup_i H_i$ is a normal crossing divisor.
Therefore the double covering
$X$ has only Gorenstein singularities and the trivial dualizing sheaf
admitting a global crepant resolution (see \cite{CM}).
A variety with this property is called a Calabi-Yau variety.

In this section, we define a double covering 
$C^\#$ of the curve $C$ and
a closed subvariety $U^{\#}$ in $C^\#\times X$
of codimension two.
A component $U$ of $U^{\#}$ defines
an algebraic correspondence between
the varieties $C$ and $X$.

Let $t$ be a point in $C$.
Then $Q_t\subset \bold P(V)$
is a cone over a smooth conic.
We define a rational curve
$\gamma_t$ in $\bold P(V^*)$ by
$$
\gamma_t=\{[W]\in \bold P(V^*)\mid
W \text{ is tangent to
the cone }Q_t
\}.
$$
If $[W]\in \gamma_t$, then $W$ contains a vertex 
$s(t)\in \bold P(V)$. Therefore the rational
curve $\gamma_t$
is contained in the hyperplane $M_{s(t)}$.
(For the notation $M_*$ and $[*]$, 
see Notation \ref{cor quad and symm}.)

\begin{proposition}
\label{tangent proposition}
The curve $\gamma_t$ is tangent to the 
hyperplanes $H_1, \dots, H_8$ for all $t\in C$.
\end{proposition}
\begin{proof}
The quadric $Q_t$ contains any of the base points $b_i$. 
Since the linear system $\Phi_L$ is regular, 
$b_i$
does not coincide with $s(t)$ for any $t\in C$.
Therefore the conic $Q_t$ is smooth at the point $b_i$ and there is
only one tangent hyperplane of $Q_t$ containing $b_i$. 
Therefore the hyperplane $H_i=M_{b_i}$ is tangent to the conic
$\gamma_i$. Thus we have proved the proposition. 
\end{proof}
\begin{corollary}
Under the assumption of Proposition \ref{tangent proposition},
the inverse image $\pi^{-1}(\gamma_t)$ is the union of two rational
curves $\tilde \gamma_t^{(1)}$ and
$\tilde \gamma_t^{(2)}$.
\end{corollary}
\begin{proof}
By Proposition \ref{tangent proposition}, 
the normalization of $\pi^{-1}(\gamma_t)$
is an etale double covering of a smooth rational curve.
Therefore it is a disjoint union of two copies
of $\gamma_t$.
\end{proof}
By attaching the irreducible components
of $\pi^{-1}(\gamma_t)$ to $t\in C$, we get 
a double covering $C^\#$ (possibly not irreducible)
of $C$.
More precisely, the covering $C^{\#}$ is defined as the Stein
factorization of the normalization of the map
$$
\{(x,t)\in X\times C\mid \pi(x)\in \gamma_t\}\to C
$$
obtained by the second projection.
The irreducible component of $\pi^{-1}(\gamma_t)$
corresponding to a point $\tilde t\in C^\#$
is written as $\tilde \gamma_{\tilde t}$.

\begin{lemma}
The covering
$C^\#$ is etale over $C$.
\end{lemma}
\begin{proof}
Since the rank of $Q_t$ is equal to $3$ for any point $t\in C$,
$\gamma_t$ is a smooth rational curve of degree two in $\bold P(V^*)$.
We will show that it is not contained in $H_i$ for any $i$.
The curve $\gamma_t$ is contained in $H_i$ if and only if the
vertex $s(t)$ of the quadric $Q_t$ coincides with one of 
a base point $b_i$
but for a smooth plane curve $C$ of degree four, any intersection of
bitangents is not contained in $C$. 
Therefore $\gamma_t$ is not contained in $H_i$.
As a consequence, $\pi^{-1}(\gamma_t)$ has
exactly two irreducible components and $C^\#$ is etale over $C$.
\end{proof}

Let $U^{\#}$ be the universal family 
$$
U^{\#}=
\{(\tilde t,x)\in C^\# \times X\mid
x\in \tilde\gamma_{\tilde t} \}.
$$
Then we have the following diagram:
$$
\begin{matrix}
U^{\#} & \subset & C^\#\times X & \overset{pr_2}\to & X \\
& & \hskip -0.2in\text{\tiny $pr_1$}\downarrow \\
& & C^\#.
\end{matrix}
$$

\begin{proposition}
\label{trivial monodromy}
The curve $C^\#$ is a union of two copies $C^{(1)}$ and $C^{(2)}$ of $C$.
\end{proposition}
\begin{proof}
To prove the proposition, it is enough to prove 
the monodromy action on
the irreducible components of 
$\pi^{-1}(\gamma_t)$ is trivial. We compute the monodromy action 
when the plane quartic curve tends to
a hyperelliptic curve. This will be done in 
Proposition \ref{trivial monodromy2}.
\end{proof}
\begin{definition}[Algebraic correspondence $\Cal U$]
\label{algebraic correspondence def}
We define $U\subset C\times X$ by the 
pull back of $U^{\#}$ by the map
$C \times X\overset{\simeq}\to C^{(1)}
\times X\subset C^{\#}\times X$.
\end{definition}

\section{Twisted cubic and hyperelliptic curves}
\label{twisted cubic section}
\subsection{Special net of quadrics}

In this subsection, we consider a special 
net of quadrics generated by
\begin{equation}
\label{special quadrics}
Q_1=
\left(
\begin{matrix}
 0& 0& 1& 0\\
 0& -2& 0& 0\\
 1& 0& 0& 0\\
 0& 0& 0& 0 
\end{matrix}\right),
Q_2=
\left(
\begin{matrix}
 0& 0& 0& -1\\
 0& 0& 1& 0\\
 0& 1& 0& 0\\
 -1& 0& 0& 0
\end{matrix}\right),
Q_3=
\left(
\begin{matrix}
 0& 0& 0& 0\\
 0& 0& 0& 1\\
 0& 0& -2& 0\\
 0& 1& 0& 0
\end{matrix}\right).
\end{equation}
Here we used the correspondence between quadratic 
polynomials and symmetric matrices
in Notation \ref{cor quad and symm}.
Then the intersection $T=Z(Q_1)\cap Z(Q_2)\cap Z(Q_3)$ is 
a twisted cubic curve defined by the image of the map
$$
s:\bold P^1 \to \bold P(V):x\mapsto (x^3:x^2:x:1).
$$
The defining equation of the singular locus 
of the net of quadrics 
$\Phi_{\langle Q_1,Q_2,Q_3\rangle}$ is equal to
$$
\det(t_1Q_1+t_2Q_2+t_3Q_3)=(t_1t_3-t_2^2)^2=0.
$$
The conic $\{t_1t_3-t_2^2=0\}$ in $\bold P(L)$ is denoted as $D$.
We choose a parameter of $D$ as $(t_1:t_2:t_3)=(1:x:x^2)$.
The singular quadratic polynomial $Q_t$ for 
$t=(t_1:t_2:t_3)=(1:x:x^2)\in D$
is written as $Q(x)$.
\begin{lemma}
For $(1:x:x^2) \in D$, we define a curve 
$\gamma_x$ 
in $\bold P(V^*)$ by
$$
\gamma_x=\{[W]\in \bold P(V^*)\mid W \text{ is tangent to }Z(Q(x))\}.
$$ 
Then $\gamma_x$
is equal to
\begin{align}
\label{equation of gamma}
\gamma_x =\{&(\xi_1:\xi_2:\xi_3:\xi_4)=
(-1:2w+x:- (w^2+2xw):w^2x)
\mid w\in \bold P^1\}.
\end{align}
Moreover, the vertex $v(x)$ of $V(Q(x))$ is 
$$(p_1:p_2:p_3:p_4)=s(x)=(x^3:x^2:x:1).$$ 
\end{lemma}
\begin{proof}
By direct computation,
we can check that the hyperplane 
\begin{align*}
W_w =\{&(p_1:p_2:p_3:p_4)\in\bold P(V)\mid \\
\nonumber
&-p_1+(2w+x)p_2 - (w^2+2xw)p_3+w^2xp_4=0
\}
\end{align*}
is tangent to the quadric $Z(Q(x))
=\{Q(x)=(p_1p_3-p_2^2)+x(p_2p_3-p_1p_4)+x^2(p_2p_4-p_3^2)=0\}$.
Since the curve $\gamma_x$ is known to be of degree two, 
we have proved the lemma.
\end{proof}

\subsection{A deformation of special net of quadrics}
Let $\Delta(u)=\{u\in \bold C\mid $ \linebreak
$\quad \mid u \mid < \epsilon\}$
be a sufficiently small disk.
The punctured dis $\Delta(u)-\{0\}$ is denoted by $\Delta^*(u)$.

\begin{proposition}
\label{generic deformation}
Let $Q_1,Q_2$ and $Q_3$ be the quadrics defined in 
(\ref{special quadrics}).

Let $\<Q_1(u),Q_2(u),Q_3(u)\>$ be a generic deformation over 
$\Delta(u)$ of 
$L=\< Q_1,Q_2,Q_3\>$ such that $Q_i(0)=Q_i$.
Then the equation of the discriminant locus 
$F(t_1,t_2,t_3,u)=\det(t_1Q_1(u)+t_2Q_2(u)+t_3Q_3(u))$
can be written as
$$
F(t_1,t_2,t_3,u)\equiv (t_1t_3-t_2^2)^2+uf(t_1,t_2,t_3)
\quad (\text{mod }u^2\bold C[t_1,t_2,t_3,u])
$$
such that $\{f(t_1,t_2,t_3)\}\cap D$ consists of distinct 
eight points
and $C':\{f(t_1,t_2,t_3)=0\}$ is a smooth plane curve.
\end{proposition}
\begin{proof}
We consider a deformation given by
\begin{align*}
& Q_1(u)=
\left(
\begin{matrix}
0 & 0 & 1 & 0 \\
0 &-2 & 0 & 0 \\
1 & 0 & 0 & 0 \\
0 & 0 & 0 & 4u
\end{matrix}
\right), \quad
 Q_2(u)=
\left(
\begin{matrix}
0 & 0 & 0 & -u-1 \\
0 & 0 & 3u+1 & 0 \\
0 & 3u+1 & 0 & 0 \\
-u-1 & 0 & 0 & 0
\end{matrix}
\right), \\
& Q_3(u)=
\left(
\begin{matrix}
4u & 0 & 0 & 0 \\
0 & 0 & 0 & 1 \\
0 & 0 &-2 & 0 \\
0 & 1 & 0 & 0
\end{matrix}
\right).
\end{align*}
Then we have $Q_i(0)=Q_i$ for $i=1,2,3$, and 
\begin{align*}
\det(t_1Q_1(u)+t_2Q_2(u)+t_3Q_3(u))\equiv
(t_1t_3-t_2^2)^2+&8(t_1^4+t_2^4+t_3^4)u \\
&\quad(\text{mod }u^2\bold C[t_1,t_2,t_3,u])
\end{align*}
The intersection of $C'=\{8(t_1^4+t_2^4+t_3^4)=0\}$ and $D$
consists of distinct eight points and $C'$ is smooth. 
Thus we have proved the proposition for a generic deformation.
\end{proof}

\subsection{Smooth family $\tilde{\Cal C}$ and its double covering 
$\Cal C^\#$}
\label{construction of family}
Let $f(t_1,t_2,t_3)$ be a generic homogeneous polynomial of 
degree four.
Then the zero locus $Z(f)$ is smooth and
the set $\{(t_1:t_2:t_3)\mid t_1t_3-t_2^2=0,f(t_1,t_2,t_3)=0\}$
consists of distinct eight points
$$
\tau_i=(1:\lambda_i:\lambda_i^2)\quad (i=1, \dots, 8).
$$
We define a family of plane curves $\Cal C$ of degree four over $\Delta(u)$ by
$$
\Cal C=\{(u,(t_1:t_2:t_3)\mid (t_1t_3-t_2^2)^2-uf(t_1,t_2,t_3)\}
$$
By Proposition \ref{generic deformation},
there is a family of nets of quadrics
$\<Q_1(u),Q_2(u),Q_3(u)\>$ in $\bold P(V)$
such that the family of associated plane quartic curve
is isomorphic to $\Cal C$.

By changing the base $\Delta(v)\to \Delta(u)$ defined by
$v\mapsto v^2=u$, we have a family $\Cal C\times_{\Delta(u)}\Delta(v)$
of plane curve in $\bold P(L)\times \Delta(v)$ over $\Delta(v)$.
Then the normalization $\tilde{\Cal C}$ of $\Cal C\times_{\Delta(u)}\Delta(v)$
is isomorphic to
$$
\{(v,(\nu:t_1:t_2:t_3))\in \Delta(v)\times 
\bold P(2,1,1,1)\mid
t_1t_3-t_2^2=\nu v, \nu^2=f(t_1,t_2,t_3)
\},
$$
where $(\nu:t_1:t_2:t_3)$ is the coordinates of the weighted projective
space $\bold P(2,1,1,1)$.
Thus we have a smooth family 
\begin{equation}
\label{general deformation}
p:\tilde {\Cal C} \to \Delta(v)
\end{equation}
of curves of genus three
over $\Delta(v)$. 
The central fiber $p^{-1}(0)=\tilde{\Cal C_0}$ is a hyperelliptic curve.
The inverse image $p^{-1}(\Delta(v)^*)$ of
$\Delta(v)^*$ is denoted as $\tilde{\Cal C^0}$.

We have a family of regular Cayley octads
$$\beta_0:\Cal B^0 \to 
\bold P(V)\times \Delta^*(v)\to
\Delta^*(v)
$$
of $\tilde{\Cal C^0} \to\Delta^*(v)$ over $\Delta^*(v)$ 
and its closure 
in $\bold P(V)\times \Delta(v)$ is denoted as
$\beta:\Cal B\to \Delta(v)$.

\begin{proposition}
\begin{enumerate}
\item
Let $v_0\neq 0$ be an element in $\Delta(v)$. 
Then the monodromy substitution on
$\Cal B^0_{v_0}=\{b_1(v_0),\dots, b_8(v_0)\}$ 
around $v=0$ is trivial. Especially $b_i(v_0)$
extends to a section $b_i(v)$ to $\bold P(V)$ over $\Delta(v)$.
\item
The point $b_i(v)$ converges to the point $s(\lambda_i)$
when $v$ tends to zero.
\end{enumerate}
\end{proposition}
\begin{proof}
(1) Choose two elements $b_i(v_0)$ and $b_j(v_0)$.
The line connecting them in $\bold P(V)\times \{v_0\}$ 
is denoted as $l_{ij}(v_0)$.
There are 28 lines of such form.
Then $l_{ij}(v_0)$ intersect the Steinerian curve at two points 
$p_{ij}^{(1)}(v_0),p_{ij}^{(2)}(v_0)$.
The intersection points correspond to
the bitangent points of the plane quartic curve
$\tilde {\Cal C_{v_0}}$ in 
$\bold P(L)$.
The divisor $p_{ij}^{(1)}(v_0)+p_{ij}^{(2)}(v_0)$ 
is one of 28 odd theta divisors.
Since the family of curves $\tilde {\Cal C}$ is smooth over $\Delta(v)$,
the monodromy action on the 28 lines $l_{ij}(v_0)$ is trivial.
Therefore the monodromy action on the set $\{b_i(v_0)\}$ is trivial.

(2) The curve $\tilde {\Cal C_0}$ is a hyperelliptic curve which
is the double covering of $D$ branched at $\tau_1, \dots, \tau_8$. 
Since the base locus of $\<Q_1(0),Q_2(0),Q_3(0)\>$ is equal to $D$,
we have $\{\tau_1,\dots, \tau_8\} \subset D$.
For a hyperelliptic curve $\tilde{\Cal C_0}$, 
an odd theta divisor is obtained
by taking sum of distinct two points in the branch points.
Therefore the family of lines $l_{ij}(v)$ connecting 
$b_i(v)$ and $b_j(v)$
converges to the line connecting
the image $s(\lambda_i)$ and $s(\lambda_j)$.
Therefore $b_i(v)$ converges to $s(\lambda_i)$.
\end{proof}

By the above proposition, the fiber $\beta^{-1}(0)$ of $\beta$ at $0$
consists of distinct 8 points 
$\{(\lambda_i^3:\lambda_i^2:\lambda_i:1)\}_{i=1,\dots,8}$.
We chose a lifting $\tilde b_i(v)=(b_{1i},\cdots,b_{4i})$ 
of the section $b_i(v)$
from $\Delta(v)$ to $V-\{0\}$
such that
$\tilde b_i(0)=(\lambda_i^3,\lambda_i^2,\lambda_i,1)$.
We consider a family of Calabi-Yau varieties $\Cal X$ 
over $\Delta(v)$ defined by
\begin{equation}
\label{CY family eq}
\Cal X_v:\eta^2=\prod_{i=1}^8
(b_{i1}(v)\xi_1+\cdots+b_{i4}(v)\xi_4)
\end{equation}
branched along the eight hyperplanes $M_{b_1(v)},\dots,M_{b_8(v)}$.
Here the equation is a weighted homogeneous equation. (The weight
of $\eta$ is four.)
We set $\Sigma=\{\tau_1,\dots, \tau_8\}\subset
p^{-1}(0)$.
For a point $t\in \tilde{\Cal C}-\Sigma$,
the inverse image $\pi^{-1}(\gamma_t)$ has exactly two irreducible components. 
These irreducible components defines an etale
double covering
$
\Cal C^{\# 0}\to \tilde{\Cal C}-\Sigma$.
By the purity of branch locus, we have an etale double
covering 
\begin{equation}
\label{family of component map}
\Cal C^{\#}\to \tilde{\Cal C}
\end{equation}
extending the above covering.

Let $\Cal U^{\# 0}\subset \Cal C^{\# 0}\times_{\Delta(v)} \Cal X$ 
be the universal family of rational curves in $\Cal X$ parameterized by 
$\Cal C^{\#0}$
and $\Cal U^\#$ the closure of $\Cal U^{\# 0}$ in
$\Cal C^{\#}\times_{\Delta(v)}\Cal X$.

\begin{proposition}
\label{trivial monodromy2}

The covering (\ref{family of component map})
is a union of two copies $\tilde{\Cal C}^{(1)}$ and 
$\tilde{\Cal C}^{(2)}$ of $\tilde{\Cal C}$
\end{proposition}
\begin{proof}
To prove the proposition,
it is enough to show
the fiber of the covering (\ref{family of component map}) at $v=0$
is a union of two copies of $\tilde{\Cal C}_0$.
The equation of the fiber of $\Cal X$ and $\tilde {\Cal C}$ at $v=0$ is 
\begin{align}
\label{CY eq}
&\Cal X_0:\eta^2=\prod_{i=1}^8
(\lambda_i^3\xi_1+\lambda_i^2\xi_2+\lambda_i\xi_3+\xi_4),
\text{ and } \\
\nonumber
&\tilde {\Cal C}_0:y^2=\prod_{i=1}^8(x-\lambda_i).
\end{align}
By restricting the equation (\ref{CY eq}) to 
the curve $\gamma_x$ given in 
(\ref{equation of gamma}),
we have
$$
\eta^2=\prod_{i=1}^8\{ (w-\lambda_i)^2(x-\lambda_i)\} .
$$
Therefore the rational curve $\gamma_x$ can be lifted to
a rational curve $\tilde\gamma_{x,y}^{(1)}$ in 
$\Cal X_0$ defined by
\begin{align}
\label{explicit eq for family of curves}
&\{(\eta:\xi_1:\xi_2:\xi_3:\xi_4)
= 
(y\prod_{i=1}^8(w-\lambda_i):
-1:2w+x: - (w^2+2xw):w^2x) \\
\nonumber
& \in \bold P(2,1,1,1) \mid 
w\in \bold P^1
\}.
\end{align}
Since the family $\{{\tilde \gamma_{x,y}}^{(1)}\}$
of rational curves 
is parameterized by $(x,y) \in\tilde{\Cal C}_0$,
we have a morphism $\tilde{\Cal C_0}\to \Cal C^{\#}_0$.
Therefore the covering $\Cal C^{\#}_0\to\tilde{\Cal C}_0$ is
a union of two copies $\tilde{\Cal C}^{(1)}_0$ and 
$\tilde{\Cal C}^{(2)}_0$
of $\tilde{\Cal C}_0$.
\end{proof}
\begin{definition}
Let $\Cal U$ be the inverse image of $\Cal U^\#$ by the open immersion 
$\tilde{\Cal C}\times \Cal X\simeq \tilde{\Cal C}^{(1)}\times \Cal X \to
\tilde{\Cal C}^{\#}\times \Cal X$.
\end{definition}

\section{Algebraic correspondence}

Let $U\subset C\times X$ be the
family of rational curves on $X$ parameterized by $C$ defined in
Definition \ref{algebraic correspondence def}.
By using the cycle class $cl(U)\in H^4(C\times X,\bold Q)(2)$,
we have the following diagram of homomorphisms:
$$
\begin{matrix}
H^1(C\times X,\bold Q) &
\overset{cl(U)}\to& H^5(C\times X,\bold Q)(2) &
\overset{pr_{2*}}\to& H^3(X,\bold Q)(1) \\
\text{\tiny $pr_1^*$}\uparrow \\
H^1(C,\bold Q)
\end{matrix}
$$
Here $V(i)$ denotes the
Tate twist of $V$.
The composite homomorphism
\begin{equation}
\label{one correspondence}
H^1(C,\bold Q) \to H^3(X,\bold Q)(1).
\end{equation}
is called the cylinder map for $U$.
We consider the family $\Cal U$ of the universal
family over $\Delta(v)$:
$$
\begin{matrix}
\Cal U & \subset & \tilde{\Cal C}\times_{\Delta(v)} \Cal X & 
\overset{pr_2}\to & \Cal X \\
& & \hskip -0.2in\text{\tiny $pr_1$}\downarrow & & \downarrow
\text{\tiny $\Xi$}\\
& & \tilde{\Cal C} &\overset{p}\to& \Delta(v) 
\end{matrix}
$$
Then we have a homomorphism $cyl$ of local systems induced by
the family of algebraic correspondences $\Cal U$:
\begin{equation}
\label{family of cylinder map}
cyl:R^1p_*\bold Q \to R^3\Xi_*\bold Q(1)
\end{equation}
\begin{theorem}
\label{correspondence thm}
The homomorphism (\ref{family of cylinder map})
is injective.
\end{theorem}
By specialization argument, we have the following corollary.
\begin{corollary}
The homomorphism
(\ref{one correspondence})
is injective.
\end{corollary}
\begin{proof}[Proof of Theorem \ref{correspondence thm}]
To prove the injectivity of the map (\ref{family of cylinder map}),
it is enough to 
show the injectivity for the fiber of (\ref{family of cylinder map})
at $v=0$.
In this case, the homomorphism $cyl$ is described by
the homomorphism
$\varphi$ defined in \S \ref{intro section}.
By the the isomorphism (\ref{basic isom}), 
we have
$$
(\tilde{\Cal C}_0\times\tilde{\Cal C}_0\times\tilde{\Cal C}_0)
/G\simeq \Cal X_0.
$$
where the finite group $G$ is defined in \S \ref{intro section}.
Let $(x,y)$ be a point in $\tilde{\Cal C}_0$.
By the equation
(\ref{explicit eq for family of curves}),
the rational curve $\tilde\gamma_{x,y}^{(1)}$ in $\Cal X_0$
defined by (\ref{explicit eq for family of curves})
is equal to the image of the map $I_{x,y}$ defined by
$$
I_{x,y}:
\tilde{\Cal C}_0
\to\Cal X_0:
(x_1,y_1)\mapsto \varphi((x_1,y_1),(x_1,-y_1),(x,y)).
$$
Therefore the universal family
$\Cal U\subset 
\tilde {\Cal C}_0\times\Cal X_0$ 
is equal to the image of the map
$$
\tilde{\Cal C}_0\times \tilde{\Cal C}_0 
\longrightarrow \tilde{\Cal C}_0\times\Cal X_0:
((x_1,y_1),(x,y))\mapsto 
((x,y),I_{x,y}(x_1,y_1))
$$
Using this description of the universal family,
the fiber $H^1(\tilde{\Cal C}_0, \bold Q) 
\to H^3(\Cal X_0,\bold Q)(1)$ of the homomorphism 
(\ref{family of cylinder map}) at $v=0$ is identified with the map
$$
H^1(\tilde{\Cal C}_0, \bold Q) 
\to \bigwedge^3H^1(\tilde{\Cal C}_0,\bold Q)(1)
\simeq H^3 (\Cal X_0,\bold Q)(1)
$$
given by $a_1\mapsto a_1\wedge \phi$,
where $\phi$ corresponds to the polarization of 
$H^1(\tilde{\Cal C_0},\bold Q)$.
Therefore it is injective.
\end{proof}
\section{Cokernel of cylinder map and infinitesimal variation of 
mixed Hodge structure}
\label{ivhs for double covering}
\subsection{Hodge structures of double coverings and Jacobian rings}

Let $M$ be the moduli space of genus three curves with level
$2N$ structures ($N\geq 2$) and $p:\Cal C\to M$ be the
universal curve over $M$.
Let $M_{he}$ be the hyperelliptic locus of $M$ and the complement
$M\backslash M_{he}$ is denoted by $M^0$.
Then the restriction $p^0:\Cal C^0=p^{-1}(M^0)\to M^0$
is a family of plane quartic curves with level $2N$ structures.
By choosing one of the family of even theta divisors
$\Theta$, we get a family 
$$
\Cal B^0 \subset 
\bold P(p^0_*\Cal O(K_{\Cal C^0/M^0}+\Theta)^*)\to M^0.
$$
of Cayley octads. The closure of $\Cal B^0$ in 
$\bold P(p_*\Cal O(K_{\Cal C/M}+\Theta)^*)$ is denoted by $\Cal B$.
Then $\Cal B$ is a disjoint union of eight copies of trivial covering
of $M$. 
By taking an etale covering $\theta:\tilde M\to M$ of $M$, 
we can choose a double covering $\tilde{\Xi}:\tilde {\Cal X}\to \tilde M$ of
$\theta^*\bold P(p_*\Cal O(K_{\Cal C/M}+\Theta))\to \tilde M$ 
branching along
the family of dual hyperplanes corresponding to $\theta^*\Cal B$. 
Moreover, by Theorem \ref{correspondence thm}, 
we have an injective morphism of local systems 
$\theta^*R^1p_*\bold Q(-1) \to R^3\tilde\Xi_*\bold Q$
by replacing $\tilde M$ by its etale covering, if necessary. 
Since the covering transformation
group $\tilde {\Cal X} \to \theta^*
\bold P(p_*\Cal O(K_{\Cal C/M}+\Theta))$
acts as $(-1)$-multiplication on $R^3\tilde {\Xi}_*\bold Q(-1)$, 
the descent data for 
$\theta^*R^1p_*\bold Q(-1)$ 
gives rise to a descent data of $\tilde {\Cal X}$.
The descended variety $\Xi:\Cal X\to M$ 
becomes a double covering of
$\bold P(p_*\Cal O(K_{\Cal C/M}+\Theta))$, 
which is called the universal
family of double coverings.

Let $\Cal C_{he}$ be the restriction of $\Cal C$ to the hyperelliptic
locus $M_{he}$.
Since the third higher direct image sheaf $R^3j_*\bold Q$
of the relative Jacobian variety $j:J(\Cal C_{he}/M_{he})\to M_{he}$
is isomorphic to $\bigwedge^3 R^1j_*\bold Q$,
the restriction of the cokernel $\Coker(cyl)$ of 
$$
cyl:R^1p_* \bold Q(-1)\to R^3\Xi_* \bold Q
$$ 
to the  hyperelliptic locus $M_{he}$ is
isomorphic to the primitive part of the third higher direct image
sheaf
$$
R^3j_{*,prim}\bold Q=\Coker(R^1j_* \bold Q(-1)\to R^3j_* \bold Q).
$$
In this section, we prove that there does not exist a polarized abelian
scheme $a:\Cal A \to M$ whose primitive part $R^3a_{*,prim}
\bold Q$ is isomorphic to $\Coker(cyl)$.

We recall computations of the infinitesimal variations of Hodge structure 
of double coverings branched along normal crossing eight hyperplanes
using Jacobian rings.
Let $\tilde b_1=(b_{11},\dots, b_{41}), \dots, 
\tilde b_8=(b_{18},\dots, b_{48})$ be non-zero elements in $V$ such that
such that the union $\cup_iH_i$ of the dual hyperplanes $H_i$
of $b_i=(b_{1i}:\cdots :b_{4i})\in\bold P(V)$ is normal crossing.
We normalize $\bold B=(^t\tilde b_1,\dots,^t \tilde b_8)$
as follows :
\begin{equation}
\label{normalized configuration}
\bold B=(b_{ij})_{ij}=\left(
\begin{matrix}
1 & 0 & 0 & 0 & 1 & 1 & 1 & 1 \\
0 & 1 & 0 & 0 & 1 & s_{11} & s_{12} & s_{13} \\
0 & 0 & 1 & 0 & 1 & s_{21} & s_{22} & s_{23} \\
0 & 0 & 0 & 1 & 1 & s_{31} & s_{32} & s_{33} \\
\end{matrix}
\right).
\end{equation}
We consider the covering $\widetilde X$ of 
$\bold P(V^*)=\{(\xi_1:\xi_2:\xi_3:\xi_4)\}$ 
in a weighted projective space
defined by
$$
\eta_i^2=b_{1i}\xi_1+\cdots + b_{4i}\xi_4. \qquad
(i=1, \dots, 8)
$$
Here the weight of $\eta_i$ is $\displaystyle \frac{1}{2}$
for $i=1, \dots,8$.
Then $\widetilde X$ is a complete intersection of Fermat type quadrics
in $\bold P(W)=\{\eta_1:\dots :\eta_8)\}$ defined by
\begin{align*}
& f_1=\eta_5^2-(\eta_1^2+\eta_2^2+\eta_3^2+\eta_4^2)=0 \\
& f_2=\eta_6^2-(\eta_1^2+s_{11}\eta_2^2+s_{21}\eta_3^2+s_{31}\eta_4^2)=0 \\
& f_3=\eta_7^2-(\eta_1^2+s_{12}\eta_2^2+s_{22}\eta_3^2+s_{32}\eta_4^2)=0 \\
& f_4=\eta_8^2-(\eta_1^2+s_{13}\eta_2^2+s_{23}\eta_3^2+s_{33}\eta_4^2)=0.
\end{align*}

We set $F(\eta,q)=\sum_{j=1}^4q_jf_j(\eta)$.
The Jacobian ideal $J(\widetilde X)$ of $\widetilde X$ is an ideal of the 
bi-graded ring $\bold C[\eta_1, \dots, \eta_8,q_1,\dots, q_4]$
generated by $\displaystyle\frac{\partial F(\eta,q)}{\partial \eta_i}$ and
$\displaystyle\frac{\partial F(\eta,q)}{\partial q_j}$. Here the bi-degrees of
$\eta_i$ and $q_i$ are $(1,0)$ and $(0,1)$, respectively.
We define the Jacobian ring $R(\widetilde X)$ of $\widetilde X$ by
the quotient ring 
$\bold C[\eta_1, \dots, \eta_8,q_1,\dots, q_4]/J(\widetilde X)$.
Then we have the following proposition (see \cite{T2}).

\begin{proposition}
\label{ivhs ci}
\begin{enumerate}
\item
The Hodge component
$H^{3-i,i}(\widetilde X)$ is identified with
$R(\widetilde X)_{2i,i}$
for $i=0,1,2,3$. The tangent space of the moduli space 
$H^1(\widetilde X, \Theta_{\widetilde X})$ is
identified with $R(\widetilde X)_{2,1}$.
\item
Under the above identification, 
the Kodaira-Spencer map 
$$
H^1(\widetilde X,\Theta_{\widetilde X})
\otimes H^{3-i,i}(\widetilde X)\to H^{2-i,i+1}(\widetilde X)
$$
is identified with 
the multiplication map of the Jacobian ring.
\item
There is an isomorphism $t:R(\widetilde X)_{6,3}\overset{\simeq}\to \bold C$
such that the cup product
$$
H^{3-i,i}(\widetilde X)\otimes
H^{i,3-i}(\widetilde X)\to 
H^{3,3}(\widetilde X)\simeq \bold C
$$ 
is identified with the following composite map
$$
R(\widetilde X)_{2i,i}\otimes
R(\widetilde X)_{2(3-i),3-i}\to
R(\widetilde X)_{6,3}\overset{t}\to \bold C.
$$
\end{enumerate}
\end{proposition}

The group $\tilde G={\mu_2}^8$ acts on the variety $\widetilde X$ and
the ring $R(\widetilde X)$ by 
$$
(\eta_i)_i\mapsto (m_i\eta_i)_i \text{ for  }(m_i)\in \tilde G.
$$
We set $G=\{(m_i)_i \in \widetilde G\mid \prod_i m_i=1\}$.
Then the double covering $X$ is isomorphic to the quotient $X=\widetilde X/G$
of $\widetilde X$ by the action of $G$.
Therefore the Hodge structure of $X$
is the invariant part of $H^3(\widetilde X)$ under the action of $G$.
\begin{proposition}
Let $R(\widetilde X)^{\tilde G}$ be the fixed part of 
$R(\tilde X)$.
Under the identification in Proposition \ref{ivhs ci}, we have
$$
H^{3-i,i}(X)\simeq R(\widetilde X)^{\tilde G}_{2i,i}
$$
for $i=0,1,2,3$. The tangent space of the moduli space 
coming from deformations of configurations of hyperplanes is
identified with $R(\widetilde X)^{\tilde G}_{2,1}$.
\end{proposition}
\begin{proof}
Let $\chi$ be the character of $\widetilde G$ defined by
$(m_i)_i \mapsto \prod_{i}m_i$ and $\Omega$
be the differential form defined by
$$
\Omega=(\sum_{i=1}^8 (-1)^i \eta_i\wedge \bigwedge_{j\neq i}d\eta_j)\wedge
(\sum_{k=1}^4 (-1)^k q_k\wedge \bigwedge_{j\neq k}dq_j).
$$
Recall that the isomorphism $R(\widetilde X)_{2i,i}\simeq H^{3-i,i}(X)$
is obtained by the composite
$$
\begin{matrix}
R(\widetilde X)_{2i,i}&\to& H^{6-i,3+i}
(\bold P^7\times \bold P^3-\{F(\eta,q)=0\})
&\simeq & H^{3-i,i}(X) \\
f(\eta,q)&\mapsto &
\frac{f(\eta,q)\Omega}{F(\eta,q)^{4+i}}
\end{matrix}
$$
Then the group $\widetilde G$ acts of on $\Omega$
and $H^3(X)$
via the character $\chi$. 
Therefore the action of $\widetilde G$ on $f$ is trivial.
\end{proof}
The invariant part of $R(\widetilde X)$ under the action of
$\tilde G$ is isomorphic to
$\bar R=\bold C[u_1, \dots, u_8, q_1,\dots, q_4]/\bar J$,
where $\eta_i^2=u_i$ and $\bar J$ is generated by
\begin{align*}
&u_5-(u_1+u_2+u_3+u_4),\quad
u_6-(u_1+s_{11}u_2+s_{21}u_3+s_{31}u_4), \\
&u_7-(u_1+s_{12}u_2+s_{22}u_3+s_{32}u_4),\quad
u_8-(u_1+s_{13}u_2+s_{23}u_3+s_{33}u_4), \\
&u_5q_1,\quad u_6q_2,\quad u_7q_3,\quad u_8q_4, \\
&u_1(q_1+q_2+q_3+q_4),\quad
u_2(q_1+s_{11}q_2+s_{12}q_3+s_{13}q_4),\\
&u_3(q_1+s_{21}q_2+s_{22}q_3+s_{23}q_4),\quad
u_4(q_1+s_{31}q_2+s_{32}q_3+s_{33}q_4).
\end{align*}
By eliminating $u_1,u_2,u_3$ and $u_4$, the ring
$\bar R$
is isomorphic to the quotient of 
$\bold C[u_1,\dots, u_4,q_1, \dots, q_4]$
by the ideal $J$ generated by
\begin{align*}
&(u_1+u_2+u_3+u_4)q_1,\quad 
(u_1+s_{11}u_2+s_{21}u_3+s_{31}u_4)q_2, \\
&(u_1+s_{12}u_2+s_{22}u_3+s_{32}u_4)q_3,
\quad (u_1+s_{13}u_2+s_{23}u_3+s_{33}u_4)q_4, \\
&u_1(q_1+q_2+q_3+q_4),\quad
u_2(q_1+s_{11}q_2+s_{12}q_3+s_{13}q_4),\\
&u_3(q_1+s_{21}q_2+s_{22}q_3+s_{23}q_4),\quad
u_4(q_1+s_{31}q_2+s_{32}q_3+s_{33}q_4),
\end{align*}
We consider a matrix $\bold B$ as in
(\ref{normalized configuration})
such that $B=\{b_1, \dots, b_8\}$ 
form a regular Cayley octad.
In this case, $s_{13},s_{23},s_{33}$ is determined 
from $s_{11},s_{21},s_{31},s_{12},s_{22},s_{32}$
as follows.
We solve the linear equation 
\begin{align*}
\left(
\begin{matrix}
1&1 &1 \\
s_{11}&s_{21}&s_{31} \\
s_{12}&s_{22}&s_{32} 
\end{matrix}
\right)
\left(
\begin{matrix}
\alpha_1 \\ \alpha_2 \\ \alpha_3
\end{matrix}
\right)
=
\left(
\begin{matrix}
1 \\ 1 \\ 1
\end{matrix}
\right),\quad
\left(
\begin{matrix}
s_{11}&s_{21}&s_{31} \\
s_{12}&s_{22}&s_{32} \\
s_{11}s_{12}&s_{21}s_{22}&s_{31}s_{32} \\
\end{matrix}
\right)
\left(
\begin{matrix}
\beta_1 \\ \beta_2 \\ \beta_3
\end{matrix}
\right)
=
\left(
\begin{matrix}
1 \\ 1 \\ 1
\end{matrix}
\right)
\end{align*}
on $\alpha_1,\alpha_2,\alpha_3,\beta_1,\beta_2,\beta_3$.
Then $s_{13},s_{23},s_{33}$ are obtained by
$$
\displaystyle s_{13}=\frac{\alpha_1}{\beta_1},\quad
s_{23}=\frac{\alpha_2}{\beta_2},\quad
s_{33}=\frac{\alpha_3}{\beta_3}.
$$
By substituting $s_{13},s_{23},s_{33}$ by the rational function of
$s_{11}, \dots, s_{32}$, the coefficient of
$F(u,p)$ is a rational function on $s_{11}, \dots, s_{32}$.
We set
\begin{align}
\label{tangent of Cayley locus}
&\tau_1=\frac{\partial F}{\partial s_{11}},\quad
\tau_2=\frac{\partial F}{\partial s_{21}},\quad
\tau_3=\frac{\partial F}{\partial s_{31}}\\
\nonumber
&\tau_4=\frac{\partial F}{\partial s_{12}},\quad
\tau_5=\frac{\partial F}{\partial s_{22}},\quad
\tau_6=\frac{\partial F}{\partial s_{32}}.
\end{align}
Let $T$ be the linear span of
the set $\{\tau_i\}_{i=1,\dots, 6}$ in $\bar R$.
Then the vector space $T$ is the tangent space of 
the moduli space of Cayley octads at $B$.

\subsection{Infinitesimal variations of Hodge structure for abelian schemes}

In this section, we use the same notations $C,X$ and $B$
of the last subsection. The point in $M$ corresponding to $C$ is denoted
by $P$.
We recall the definition of infinitesimal variations of polarized Hodge structure.
\begin{definition}[Infinitesimal variation of 
polarized Hodge structure \cite{CG}]
A triple $(\oplus H^{i,j},\theta,(\ ,\ ))$ consisting of
\begin{enumerate}
\item 
a vector space $\oplus_{i+j=q}H^{ij}$ with a linear map
$\theta:H^{i,j}\to H^{i-1,j+1}$, and
\item
a system of $(-1)^q$-symmetric perfect pairing 
$(\ ,\ ):H^{i,j}\otimes H^{j,i}\to \bold C$
\end{enumerate}
is called an infinitesimal variation of polarized Hodge structure of weight $q$ if
$(\theta(a),b)+(a,\theta(b))=0$.
\end{definition}

To a variation of Hodge structures $\Cal H$ of weight $q$ 
over $M$ and a tangent
vector at a point $P$ in $M$, we can naturally associate
an infinitesimal variation of polarized Hodge structure
on the direct product $\oplus_{i+j=q} \Cal H_{P}^{i,j}$
of Hodge components of the fiber $\Cal H_P$.
The variation of Hodge structures $R^3\Xi_*\bold Q$ 
and a tangent vector $\theta$ at $P$
give rise to the following infinitesimal variation of 
polarized Hodge structure.
$$
\theta_{X}:H^{3-i,i}(X)\to H^{2-i,i+1}(X).
$$

Similarly, the variation of Hodge structures $R^1p_*\bold Q$ gives rise to
the following infinitesimal variation of polarized Hodge 
structure.
$$
\theta_{C}:H^{1,0}(C)\to H^{0,1}(C).
$$
Since the homomorphism 
(\ref{family of cylinder map})
comes form a family of algebraic correspondences,
the two infinitesimal variations $\theta_{C}, \theta_{X}$ are
compatible. Therefore we have the following commutative diagram.
\begin{align}
\label{compatibility of ivhs}
\begin{matrix}
H^{1,0}(C)& \overset{\theta_{C}}\to & H^{0,1}(C) \\
cyl\downarrow & & \downarrow cyl\\
H^{2,1}(X)& \underset{\theta_{X}}\to & H^{1,2}(X) 
\end{matrix}
\end{align}
The vertical map $cyl$ is injective by Theorem \ref{correspondence thm}. The cokernel
of the map $cyl$ is 
an infinitesimal variation of polarized Hodge structure,
which is denoted by $(H^3_{prim}(X),\theta_X)$.

Let $H^{2,1}_{prim}(X)$ and $H^{1,2}_{prim}(X)$ be the orthogonal complements
of the image of $H^{0,1}(C)$ and $H^{1,0}(C)$, respectively.
Then the restriction of the map $\theta_{X}$ to
$H^{3,0}\oplus H^{2,1}_{prim}\oplus H^{1,2}_{prim}\oplus H^{0,3}$
is an infinitesimal variation of polarized Hodge structure and
it isomorphic to $(H^3_{prim}(X),\theta_X)$.
Then the map
$$
TM_P \to Hom(H^{3,0}(X),H^{1,2}_{prim}(X))
$$
defined by $\theta\mapsto \theta_X\circ\theta_X$
is a 
$Hom(H^{3,0}(X_p),H^{1,2}_{prim}(X))$-valued quadratic form on $TM_P$,
which is denoted as $\Cal Q_X$.

We consider infinitesimal variations of Hodge structure arising
from weight one variations.
Let $H^{1,0},H^{0,1}$ be 3-dimensional $\bold C$-vector spaces
and $(,)$ be a non-degenerate pairing $H^{1,0}\otimes H^{0,1} 
\to \bold C$.
Using the isomorphism $(H^{1,0})^*\overset{\iota}\to H^{0,1}$ and
$(H^{0,1})^*\overset{\iota}\to H^{1,0}$ obtained by this pairing,
we have the following isomorphism :
$$
Hom(H^{0,1}, H^{1,0})=Hom((H^{1,0})^*, (H^{0,1})^*)\overset{\iota}\to Hom(H^{0,1}, H^{1,0}).
$$
The image of $\theta$ under the above isomorphism is denoted by $\ ^t\theta$.
For a linear map $\theta\in Hom(H^{1,0}, H^{0,1})$, the triple
$(\bold H,\theta)=(H^{1,0}\oplus H^{0,1},\theta:H^{1,0}\to H^{0,1},(,))$ 
is an infinitesimal variation of Hodge
structure if and only if $\theta=\ ^t\theta$.
We set
$$
T=\{\theta \in Hom(H^{1,0}, H^{0,1}) \mid \theta=\ ^t\theta \}.
$$
The third exterior product $(\bigwedge^3\bold H,\theta^{(3)})$ 
$$
\wedge^3H^{1,0}\overset{\theta^{(3)}}\to \wedge^2H^{1,0}\otimes H^{0,1} 
\overset{\theta^{(3)}}\to
H^{1,0}\otimes\wedge^2H^{0,1} 
\overset{\theta^{(3)}}\to \wedge^3H^{0,1}
$$
of $(\bold H,\theta)$ is given 
by the formula 
$$
\theta^{(3)}(a_1\wedge a_2\wedge a_3)=
\theta(a_1)\wedge a_2\wedge a_3+a_1\wedge \theta(a_2)\wedge a_3+
a_1\wedge a_2 \wedge \theta(a_3)
$$
for $a_1,a_2,a_3\in H^{0,1}\oplus H^{1,0}$.
Let $\phi$ be the element in $\bigwedge^2 \bold H$ corresponding to the
identity map of $Hom(H^{0,1},H^{0,1})$
via the following isomorphism:
$$
Hom(H^{0,1},H^{0,1})\simeq (H^{0,1})^*\otimes H^{0,1}
\simeq H^{1,0}\otimes H^{0,1} 
\subset \bigwedge^2
\bold H.
$$
We set $H^{2,1}_{prim}$ and $H^{1,2}_{prim}$ the orthogonal complements
of $\phi\wedge H^{0.1}\subset H^{1,2}$ and 
$\phi\wedge H^{1.0}\subset H^{2,1}$, respectively.  
Then we have an infinitesimal variation of Hodge structure
$$
(\bigwedge^3\bold H)_{prim}=
(\bigwedge^3 H^{1,0}\oplus H^{2,1}_{prim}\oplus H^{1,2}_{prim}\oplus
\bigwedge^3 H^{0,1}, \theta^{(3)},(,)).
$$
which is isomorphic to 
$\Coker(\bold H\to \bigwedge^3 \bold H)$.
Then the map
$$
\Cal Q:T \to Hom(\bigwedge^3 H^{1,0},H^{1,2}_{prim}):\theta\mapsto
\theta^{(3)}\circ \theta^{(3)}
$$
is a $Hom(\bigwedge^3 H^{1,0},H^{1,2}_{prim})$-valued quadratic form on $T$.
We have the following proposition.

\begin{proposition}
\label{condition for exterior product}
There exists an non-zero element $\theta\in T$ such that $\Cal Q(\theta)=0$.
\end{proposition}
\begin{proof}
Let $e_1,e_2,e_3$ be a basis of $H^{1,0}$ and $f_1,f_2,f_3$ be
the dual base of $H^{0,1}$.
Let $\theta\in Hom(H^{1,0},H^{0,1})$ be the map defined by
$\theta(e_1)=f_1,\theta(e_2)=\theta(e_3)=0$.
Then $\theta \in T$ and we have
$$
(\theta^{(3)}\circ \theta^{(3)})(e_1\wedge e_2 \wedge e_3)=
\theta^{(3)}(f_1\wedge e_2 \wedge e_3)=0.
$$
\end{proof}

\begin{proposition}
\label{existence of non exterior ivhs}
There exists a regular Cayley octad $B=\{b_1, \dots,b_8\}$ with the
following property. 
\begin{enumerate}
\item
The value $\Cal Q(\theta)_X$ of $\Cal Q_X$ at 
$\theta$ is non-zero for all $\theta \neq 0$.
\item
The map $TM_P\to Hom(H^{3,0}(X),H^{2,1}(X))$ is an isomorphism.
\end{enumerate}
\end{proposition}
We give an explicit example of Cayley octad 
with the above properties in the next section. 
By Proposition \ref{condition for exterior product}
together with Proposition 
\ref{existence of non exterior ivhs},
we have the following theorem.
\begin{theorem}
\label{revised second main theorem}
There exists no polarized abelian scheme $a:\Cal A \to M$ of relative
 dimension three, whose primitive part
$R^3a_{prim,*}\bold Q$ of $R^3a_*\bold Q$ is isomorphic to
$\Coker(cyl)$
as variations of Hodge structures.
\end{theorem}
\begin{proof}
Assume that there exists an abelian scheme $a:\Cal A \to M$
such that $R^3a_{*,prim}\bold Q \simeq \Coker(cyl)$.
Let $P'\in M_{et}$ be a point corresponding to a hyperelliptic curve.
Then the image of the fundamental group $\pi_1(M,P')$ under the monodromy
representation on
$Aut((R^3a_{*,prim}\bold Q)_{P'})$ contains the image of the
fundamental group $\pi_1(M_{he},P')$ of $M_{he}$, which is finite index
in $Sp(3,\bold Z)$ by \cite{TS}. Therefore there is at most one polarization
on $R^3a_{*,prim}\bold Q$ and $\Coker(cyl)$ and 
the isomorphism of variations of Hodge structures between
$R^3a_{prim,*}\bold Q$ and $\Coker(cyl)$
induces an isomorphisms of infinitesimal variation of 
polarized Hodge structure for any $\theta\in TM_P$ at any $P\in M$.

We consider the infinitesimal variation of 
polarized Hodge structure at the point
$P$ corresponding to the Cayley octad satisfying the condition of
Proposition \ref{existence of non exterior ivhs}.
By the condition (2) and
Proposition \ref{condition for exterior product}, 
there exists non-zero tangent vector $\theta\in MT_P$
such that $\theta^{(3)}\circ \theta^{(3)}=0$ in
$Hom(H^{3,0}(\Cal A_P),H^{2,1}_{prim}(\Cal A_P))$.
But this implies
$\theta_X\circ \theta_X=0$, 
which contradicts to
the condition (1) of Proposition \ref{existence of non exterior ivhs}.
\end{proof}
\section{Example satisfying the condition in Proposition
\ref{existence of non exterior ivhs}}
\label{computer calculation}
We compute examples of the Jacobian ring using Gr\"obner base 
under the computer program Maple.
We show that there exists a Cayley octad satisfying the condition
in Proposition \ref{existence of non exterior ivhs}.
Under the notation of the last section, we set
$$
s_{11}=-1,s_{21}=3,s_{31}=4,s_{12}=-3,s_{22}=2,s_{32}=3.
$$
Then $s_{13}=\frac{33}{5},s_{23}=22,s_{33}=99$.

We identify the space $H^{2,1}(X_p)$ with the degree $(2,2)$ part 
of \linebreak
$\bold Q[u_1, \dots, u_4,q_1,\dots, q_4]/J$.
Using graded reverse lexicographic
order of \linebreak
$u_4,u_3,u_2,u_1,q_4,q_3,q_2,q_1$, 
$H^{2,1}(X_p)$ is isomorphic to the vector spanned by the class of
\begin{align*}
&M_1=u_2^2p_2^2,\quad M_2=u_2u_1p_2^2,\quad M_3=u_1^2p_2^2,\quad M_4=u_3^2p_1^2, 
\quad M_5=u_3u_2p_1^2, 
\quad \\
&M_6=u_2^2p_1^2,\quad M_7=u_3u_1p_1^2,\quad M_8=u_2u_1p_1^2,
\quad M_9=u_1^2p_1^2.
\end{align*}
Therefore we have an isomorphism
$r:H^{2,1}(X_p)\overset{\simeq}\to 
\bold C M_1\oplus\cdots \oplus \bold C M_9.$
By a simple calculation, 
the tangent vectors
$\tau_1, \dots, \tau_6$ 
in (\ref{tangent of Cayley locus})
can be computed as
$$
\tau_1=
{q_{2}}\,{u_{2}} + {q_{4}}\,( - 6\,{u_{2}} - {\displaystyle 
\frac {73}{2}} \,{u_{3}} - 420\,{u_{4}}),
$$
and so on. Let $\theta =w_1\tau_1+ \cdots +w_6 \tau_6$
be a linear combination of $\tau_1, \cdots \tau_6$.
Then $r(\theta^2)$ is a $\bold C^6$-valued quadratic form
on $w_1,\dots, w_6$.
The coefficients $f_1, \dots, f_9$ of $M_1, \dots, M_9$ 
are computed as
\begin{align*}
f_1 = & 
- {\displaystyle \frac {24757}{1240}} \,w_5^{2} 
- {\displaystyle \frac {1197}{155}} \,w_5\,w_6 
+ {\displaystyle \frac {399}{62}} \,w_2\,w_6 
+ {\displaystyle \frac {247}{124}} \,w_2\,w_4 
+ {\displaystyle \frac {4351}{310}} \,w_2\,w_5 \\
& 
+ {\displaystyle \frac {3857}{1240}} \,w_3^{2}  
- {\displaystyle \frac {1463}{1488}} \,w_3\,w_4 
+ {\displaystyle \frac {399}{62}} \,w_3\,w_5 
- {\displaystyle \frac {31787}{3720}} \,w_3\,w_6 
+ {\displaystyle \frac {29735}{35712}} \,w_4^{2} \\
&
- {\displaystyle \frac {209}{31}} \,w_4\,w_5  
+ {\displaystyle \frac {15295}{8928}} \,w_4\,w_6 
+ {\displaystyle \frac {133}{248}} \,w_1\,w_3 
+ {\displaystyle \frac {95}{62}} \,w_1\,w_2 
+ {\displaystyle \frac {12255}{992}} \,w_1^{2} \\
&
- {\displaystyle \frac {20653}{2976}} \,w_1\,w_4 
- {\displaystyle \frac {95}{62}} \,w_2^{2}   
- {\displaystyle \frac {798}{155}} \,w_2\,w_3 
+ {\displaystyle \frac {52535}{8928}} \,w_6^{2} 
- {\displaystyle \frac {2261}{1488}} \,w_1\,w_6 \\
&
+ {\displaystyle \frac {323}{62}} \,w_1\,w_5,
\end{align*}
and so on. Then by using Gr\"obner basis, we conclude
that 
$$
\{(w_1, \dots, w_6)\mid  f_1(w_1, \dots, w_6)=\dots =
f_9(w_1, \dots,w_6)=0\}=
\{(0,\dots, 0)\}.
$$
Thus we have a regular Cayley octad $B$ satisfying the conditions of
Proposition \ref{existence of non exterior ivhs}.

\end{document}